\documentclass[12pt,a4paper]{article}

\usepackage{lipsum}
\usepackage[a4paper,margin=1in]{geometry}
\usepackage{lineno}
\usepackage[english]{babel}
\usepackage{amsmath,amsthm,amsfonts,amssymb,amscd,cite}
\usepackage[colorlinks=true, linkcolor=red, filecolor=brown, citecolor=green, urlcolor=blue,bookmarksnumbered, plainpages, backref]{hyperref}

\usepackage{graphicx}
\usepackage{float}
\def\E{\mathcal{E}}

\def\Tr{\mathrm{Tr}}
\def\diag{\mathrm{diag}}

 \newtheorem{theorem}{Theorem}
 \newtheorem{lemma}{Lemma}
 \newtheorem{corollary}{Corollary}
 
\title{Vertex weighted Laplacian graph energy and other topological indices}
\author{Reza Sharafdini,\quad\quad H. Panahbar\\
Department of Mathematics, Faculty of Science,\\ Persian Gulf University,
Bushehr 7516913817, Iran\\ [3mm]
e-mail: sharafdini@pgu.ac.ir\\ [1mm]
}
\date{}
\begin{document}
\maketitle


\begin{abstract}
  Let $G$ be a graph with a vertex weight $\omega$ and the vertices $v_1,\ldots,v_n$. The Laplacian matrix of $G$ with respect to $\omega$ is defined as
  $L_\omega(G)=\diag(\omega(v_1),\cdots,\omega(v_n))-A(G)$, where $A(G)$ is the adjacency matrix of $G$. Let $\mu_1,\cdots,\mu_n$ be eigenvalues of $L_\omega(G)$. Then the Laplacian energy of $G$ with respect to $\omega$ defined as $LE_\omega (G)=\sum_{i=1}^n\big|\mu_i - \overline{\omega}\big|$, where  $\overline{\omega}$ is the average of $\omega$, i.e., $\overline{\omega}=\dfrac{\sum_{i=1}^{n}\omega(v_i)}{n}$. In this paper we consider several natural vertex weights of $G$ and obtain some inequalities between the ordinary 
  and Laplacian energies of $G$ with corresponding vertex weights. Finally, we apply our results to the molecular graph of toroidal fullerenes (or achiral polyhex nanotorus).\\[5mm]
\end{abstract}

\section{Introduction}

In this paper, we are concerned with simple graphs. Let $G$ be a simple graph, with nonempty vertex set $V(G)= \{v_1,\ldots,v_n\}$ and edge set $E= \{e_1,\ldots,e_m\}$. That is to say, $G$ is a simple $(n,m)$-graph.
Let $\omega$ be a vertex weight of $G$, i.e., $\omega$ is a function from $V(G)$ to the set of positive real numbers. In this case, we say that $G$ is a graph with a vertex weight $\omega$. A vertex weight $\omega$ could be a constant function. In this case, we say $G$ is $\omega$-regular. Namely, $G$ is $\omega$-regular if for any $u,v\in V(G)$, $\omega(u)=\omega(v)$.
Observe that a well-known vertex weight of a graph is the vertex degree assigning to each vertex its degree. Let us denote it by $deg$.

The diagonal matrix of order $n$ whose $(i,i)$-entry is $\omega(v_i)$ is called the diagonal vertex weight matrix of $G$ with respect to  $\omega$   and is denoted by $D_\omega(G)$, i.e., $D_\omega(G)=\mathrm{diag}(\omega(v_1),\ldots,\omega(v_n))$ . The adjacency matrix $A(G)= (a_{ij})$ of $G$ is a $(0,1)$-matrix defined by $a_{ij}= 1$ if and only if the vertices $v_i$ and $v_j$ are adjacent. Then the matrices $L_{deg}(G)=D_{deg}(G)-A(G)$
and $L^\dag_{deg}(G)=A(G)+D_{deg}(G)$ are called Laplacian  and  signless Laplacian  matrix of $G$, respectively (see \cite{Gro-Mer1994}, \cite{Gro-Mer-sun1990}, \cite{Mer1995}, \cite{Mer1994}, \cite{Mohar1991} and \cite{Mohar2004}). These matrices was generalized  for arbitrary vertex weighted graphs (see \cite{Sh-Pa-weiLa} and \cite{Sh-Ataei-La}). Let $G$ be a simple graph with the vertex weight  $\omega$. Then we shall call the matrices $L_\omega(G)=D_\omega(G)-A(G)$
and $L^\dag _\omega(G)=A(G)+D_\omega(G)$  the weighted  Laplacian  and the weighted signless   Laplacian matrix of $G$ with respect to the vertex weight $\omega$.

Let  $X = \{x_1 ,x_2 ,...,x_n\}$ be a data set of real numbers. The \emph{mean absolute
deviation} (often called the mean deviation) $\mathrm{MD}(X)$ and variance $\mathrm{Var}(X)$ of $X$ is defined as
\begin{linenomath*}
\begin{equation*}
\mathrm{MD}(X)= \dfrac{1}{n}\sum_{i=1}^n|x_i-\overline{x}|,\quad\quad \mathrm{Var}(X)= \dfrac{1}{n}\sum_{i=1}^n(x_i-\overline{x})^2,
\end{equation*}
\end{linenomath*}
where $\overline{x}=\dfrac{\sum_{i=1}^{n}x_i}{n}$ is the arithmetic mean of the distribution.
Note that an easy application of the Cauchy-Schwarz inequality gives that the mean deviation is a
lower bound on the standard deviation (see \cite{CaversThesis}).
\begin{linenomath*} \begin{equation}\label{MDVAR}
  \mathrm{MD}(X)\leq \sqrt{\mathrm{Var}(X)}.
\end{equation}\end{linenomath*}
The mean deviation and variance of $G$ with respect to $\omega$, denoted by $\mathrm{MD}_\omega(G)$ and  $\mathrm{Var}_\omega(G)$, respectively, is defined as
\begin{linenomath*} \begin{equation*}
\mathrm{MD}_\omega(G)=\mathrm{MD}(\omega(v_1),\ldots,\omega(v_n)),\quad \quad \mathrm{Var}_\omega(G)=\mathrm{Var}(\omega(v_1),\ldots,\omega(v_n)).
\end{equation*} \end{linenomath*}
It follows from  Eq. \eqref{MDVAR} that $\mathrm{MD}_\omega(G)\leq \sqrt{\mathrm{Var}_\omega(G)}$.
It is worth mentioning that $\mathrm{Var}_{\deg}(G)$ is well-investigated graph invariant (see \cite{Bell} and \cite{Gut-Paule2002}).
Let $\lambda_1, \lambda_2,\ldots , \lambda_n$ be eigenvalues of the adjacency matrix $A(G)$ of graph $G$. It is  known that $\sum_{i=1}^{n}\lambda_i=0$. The notion of the energy $\E(G)$ of an $(n,m)$-graph $G$ was introduced by Gutman in connection with the $\pi$-molecular energy (see \cite{I.Gutman}, \cite{Gutman}, \cite{Gutm-zar} and \cite{Ind-Vijay}). It is defined as
\begin{linenomath*} \begin{equation*}\E(G)=\sum_{i=1}^n \vert\lambda_i\vert=n \mathrm{MD}(\lambda_1, \lambda_2,\ldots , \lambda_n) .\end{equation*} \end{linenomath*}
For details of the theory of graph energy see \cite{Gutman}, \cite{Gut-Pol} and \cite{LiShiGut}.


Let $n\ge\mu_1, \mu_2,\ldots , \mu_n=0$ be eigenvalues of Laplacian matrix $L_{deg}(G)$ of graph $G$.
It is  known that $\sum_{i=1}^{n}\mu_i=2m$.
Gutman and Zhou defined the Laplacian energy of an $(n,m)$-graph $G$ for the first time (see \cite{Gut-Zhou-Lap} ) as
\begin{linenomath*} \begin{equation*} LE(G)=\sum_{i=1}^n\Big\vert\mu_i-\frac{2m}{n}\Big\vert=n\mathrm{MD}(\mu_1,\ldots, \mu_n).\end{equation*} \end{linenomath*}
Numerous results on the Laplacian energy have been obtained, see for instance
\cite{T.Ale}, \cite{Dasa.Moj.Gut}, \cite{N.N.M}, \cite{I.N.M}, \cite{Romatch2009}, \cite{KeyFanLAA2010} and \cite{B.I}.
Note that in the definition of Laplacian energy $\dfrac{2m}{n}$ is the average vertex degree of $G$. This motivates us to extend their definition to the graphs equipped with an arbitrary vertex weight. Let $G$ be
 a graph with the vertex set $V= \{v_1,\ldots,v_n\}$ and with an arbitrary vertex weight $\omega$.
Let $\mu_1, \mu_2,\ldots , \mu_n$ be eigenvalues of the weighted Laplacian matrix $L_\omega(G)$ of graph $G$ with respect to the vertex weight $\omega$. Then we  \cite{Sh-Pa-weiLa} proposed the Laplacian energy $LE_\omega(G)$ of $G$ with respect to the vertex weight $\omega$ as
\begin{linenomath*}
\begin{equation}\label{eqn:wLE}
 LE_\omega (G)=\sum_{i=1}^n\big|\mu_i - \overline{\omega}\big|=n\mathrm{MD}(\mu_1,\ldots, \mu_n),
\end{equation}
\end{linenomath*}
where
\begin{linenomath*}
\begin{equation*}
\overline{\omega}=\dfrac{\sum_{i=1}^{n}\omega(v_i)}{n} \quad \mbox{and} \quad \sum_{i=1}^{n}\mu_i=n\overline{\omega}.
\end{equation*}
\end{linenomath*}
Note that $LE_{deg}(G)=LE(G)$.

Let $G$ be a graph with an arbitrary vertex weight $\omega$. Some inequalities between $\E(G)$ and $LE_\omega(G)$  were established in \cite{Sh-Ataei-La}; and therein, the following three theorems were proved.
\begin{theorem}\label{thm:MD}
Let $G$ be a connected $(n,m)$-graph with a vertex weight $\omega$. Then
\begin{equation}\label{eqn:Lap}
 LE_\omega(G)\leq n\mathrm{MD}_\omega(G)+\E(G).
\end{equation}
Moreover the equality in \eqref{eqn:Lap} holds if and if $G$ is  $\omega$-regular.
\end{theorem}

\begin{theorem}\label{th:Reg}
Let $G$ be a bipartite graph with a vertex weight $\omega$. Then
\begin{equation}\label{eqn:Reg-bip}
LE_\omega(G)\ge \E(G).
\end{equation}
Moreover, the equality in \eqref{eqn:Reg-bip} holds if and only if $G$ is a  $\omega$-regular graph.
\end{theorem}

\begin{theorem}\label{th:}
Let $G$ be a bipartite (m,n)-graph with a vertex weight $\omega$. Then
\begin{equation}\label{eqn:}
\max\Big\{n\mathrm{MD}_\omega(G), \E(G)\Big\}\leq LE_\omega(G)\leq n\mathrm{MD}_\omega(G)+\E(G).
\end{equation}
\end{theorem}
In this paper we aim to apply the above theorems to graphs with some natural vertex weights and establish relationships between some graph invariants and Laplacian graph energy with respect to corresponding vertex weight.

\section{Main Results}

Having a molecule, if we represent atoms by vertices and bonds by edges, we obtain a
molecular graph. Graph theoretic invariants of molecular graphs, which predict
properties of the corresponding molecule, are known as topological indices.
The oldest topological index is the Wiener index, which was introduced in 1947. Since then several topological indices have been proposed to predict characteristics of chemical compounds, like physio-chemical, pharmacologic, toxicological and other biological properties.
In this article we deals with Wiener index, total eccentricity index and first Zagreb index.

\subsection{Wiener index}
Let $G$ be a connected graph.
Given two vertices $u$ and $v$ in $V(G)$, the distance between $u$ and $v$, denoted by $d ( u , v ) = d_G ( u , v )$
is the length of a shortest path connecting them. The Wiener index $W ( G )$ of a connected graph $G$ is defined
to be the sum of distances between any two unordered pair of vertices of $G$, i.e.,

\[W(G)=\sum_{\{u,v\}\subseteq V(G)}d_G ( u , v )=\frac{1}{2}\sum_{u,v\in V(G)}d_G ( u , v ).\]
The transmission $\Tr(v)$ of a vertex $v$ is defined to be the sum of the distances from $v$ to all other vertices in $G$, i.e.,
\[\Tr(v)=\sum_{u\in V(G)}d_G (u , v).\]
It is clear that
\[W(G)=\frac{1}{2}\sum_{v\in V(G)}\Tr(v).\]
A connected graph $G$ is said to be $k$-transmission regular if $\Tr ( v ) = k$ for every vertex $v\in V(G)$.
The transmission regular graphs are exactly the distance-balanced graphs introduced in \cite{Handa}. They are
also called self-median graphs \cite{Cabello}. We may consider the transmission of an arbitrary vertex as a vertex weight  with the average    $\overline{\Tr}=\frac{2W(G)}{n}$. In this point of view,  it follows from \eqref{eqn:wLE} that
\begin{linenomath*}
\begin{equation}\label{eqn:TrLE}
 LE_{\Tr} (G)=\sum_{i=1}^n\Big|\mu_i - \frac{2W(G)}{n}\Big|=n\mathrm{MD}(\mu_1,\ldots, \mu_n).
\end{equation}
\end{linenomath*}



\begin{theorem}\label{thm:MD-Lap-Tr}
Let $G$ be a connected graph with $n$ vertices.  Then
\begin{equation}\label{eqn:MD-Lap-Tr}
 LE_\Tr(G)\leq n\mathrm{MD}_\Tr(G)+\E(G).
\end{equation}
Moreover the equality in \eqref{eqn:MD-Lap-Tr} holds if and if $G$ is  transmission regular.
\end{theorem}

\begin{theorem}\label{thm:MD-Lap-Tr}
Let $G$ be a bipartite graph. Then
\begin{equation}\label{eqn:MD-Lap-Tr-bi}
LE_{\Tr}(G)\ge \E(G).
\end{equation}
Moreover, the equality in \eqref{eqn:MD-Lap-Tr-bi} holds if and only if $G$ is transmission regular.
\end{theorem}

\begin{theorem}\label{thm:MD-Lap-Tr-bi-low}
Let $G$ be a bipartite graph with $n$ vertices.  Then
\begin{equation}\label{eqn:MD-Lap-Tr-bi-low}
\max\Big\{n\mathrm{MD}_{\Tr}(G), \E(G)\Big\}\leq LE_{\Tr}(G)\leq n\mathrm{MD}_{\Tr}(G)+\E(G).
\end{equation}
\end{theorem}

\subsection{Zagreb indices}
First Zagreb index of a graph $G$ is defined as
\[M_1(G)=\sum_{u\in V(G)}deg_G(u)^2\]
For a graph $G$, denote by $t(u)$ the 2-degree of vertex $u$, which is the sum of the degrees of
the vertices adjacent to $u$;
A graph is said to be 2-degree regular if  $t(u)$ is constant for each $u$. It is known that \[M_1(G)=\sum_{u\in V(G)}t(u),\]


We may consider the 2-degree of an arbitrary vertex as a vertex weight  with the average    $\overline{t}=\frac{M_1(G)}{n}$. In this point of view,  it follows from \eqref{eqn:wLE} that

\begin{linenomath*}
\begin{equation}
 LE_{t} (G)=\sum_{i=1}^n\big|\mu_i - \frac{M_1(G)}{n}\big|=n\mathrm{MD}(\mu_1,\ldots, \mu_n).
\end{equation}
\end{linenomath*}

\begin{theorem}\label{thm:MD-Lap-1Zagreb}
Let $G$ be a connected graph with $n$ vertices.  Then
\begin{equation}\label{eqn:MD-Lap-1Zagreb}
 LE_t(G)\leq n\mathrm{MD}_t (G)+\E(G).
\end{equation}
Moreover the equality in \eqref{eqn:MD-Lap-1Zagreb} holds if and if $G$ is 2-degree regular.
\end{theorem}

\begin{theorem}\label{thm:MD-Lap-1Zagreb-bi}
Let $G$ be a bipartite graph with $n$ vertices.  Then
\begin{equation}\label{eqn:MD-Lap-1Zagreb-bi}
LE_t (G)\ge \E(G).
\end{equation}
Moreover, the equality in \eqref{eqn:MD-Lap-1Zagreb-bi} holds if and only if $G$ is 2-degree regular.
\end{theorem}

\begin{theorem}\label{thm:MD-Lap-1Zagreb-bi-low}
Let $G$ be a bipartite graph with $n$ vertices.  Then
\begin{equation}\label{eqn:MD-Lap-1Zagreb-bi-low}
\max\Big\{n\mathrm{MD}_t (G), \E(G)\Big\}\leq LE_t (G)\leq n\mathrm{MD}_t (G)+\E(G).
\end{equation}
\end{theorem}

Let us define $deg^2(u)=deg_G(u)^2$, the square vertex degree of $u$.
So we may consider $deg^2$ as a vertex weight of $G$ with the average $\overline{deg^2}=\frac{M_1(G)}{n}$. From this point of view, a graph is square vertex degree regular if and only if it is vertex degree regular. In this point of view,  it follows from \eqref{eqn:wLE} that

\begin{linenomath*}
\begin{equation}
 LE_{deg^2} (G)=\sum_{i=1}^n\big|\mu_i - \frac{M_1(G)}{n}\big|=n\mathrm{MD}(\mu_1,\ldots, \mu_n).
\end{equation}
\end{linenomath*}

\begin{theorem}\label{thm:MD-1Zagreb1}
Let $G$ be a connected graph with $n$ vertices.  Then
\begin{equation}\label{eqn:MD-1Zagreb1}
 LE_{deg^2}(G)\leq n\mathrm{MD}_{deg^2} (G)+\E(G).
\end{equation}
Moreover the equality in \eqref{eqn:MD-1Zagreb1} holds if and if $G$ is vertex degree regular.
\end{theorem}

\begin{theorem}\label{th:MD-1Zagreb1-bi}
Let $G$ be a bipartite graph. Then
\begin{equation}\label{eqn:MD-1Zagreb1-bi}
LE_{deg^2} (G)\ge \E(G).
\end{equation}
Moreover, the equality in \eqref{eqn:MD-1Zagreb1-bi} holds if and only if $G$ is vertex degree regular.
\end{theorem}

\begin{theorem}\label{thm:MD-1Zagreb1-bi-low}
Let $G$ be a bipartite graph with $n$ vertices.  Then
\begin{equation}\label{eqn:MD-1Zagreb1-bi-low}
\max\Big\{n\mathrm{MD}_{deg^2} (G), \E(G)\Big\}\leq LE_{deg^2} (G)\leq n\mathrm{MD}_{deg^2} (G)+\E(G).
\end{equation}
\end{theorem}

\subsection{Total eccentricity index}
The eccentricity $\varepsilon (u)$ of the vertex $u$ of a connected graph $G$ is the distance from $u$ to any vertex farthest away from it in $G$, i.e., $\varepsilon (u)=\max_{v\in V(G)} d(u,v)$. The maximum eccentricity over all vertices of $G$ is called the diameter of $G$ and is denoted by $D(G)$; the minimum
eccentricity among the vertices of $G$ is called the radius of $G$ and is denoted by $R(G)$.
The set of all vertices of minimum eccentricity is called the center of $G$. A connected graph $G$ is called self-centred if $\varepsilon (u)=R(G)$  for each $u\in V(G)$. The total eccentricity index of a connected graph  $G$, denoted by $\zeta(G)$,  is defined as the sum of eccentricities of vertices of $G$, i.e., $\zeta(G)=\sum_{u\in V(G}\varepsilon (u).$

One may consider the eccentricity of a vertex as a vertex weight of $G$ with the average $\overline{\varepsilon}=\frac{\zeta(G)}{n}$. From this point of view, a graph is $\varepsilon$-regular if and only if it is self-centred. In this point of view,  it follows from \eqref{eqn:wLE} that

\begin{linenomath*}
\begin{equation}
 LE_{\varepsilon}(G)=\sum_{i=1}^n\big|\mu_i - \frac{\zeta(G)}{n}\big|=n\mathrm{MD}(\mu_1,\ldots, \mu_n).
\end{equation}
\end{linenomath*}

\begin{theorem}\label{thm:MD-Total-ecc}
Let $G$ be a connected graph with $n$ vertices.  Then
\begin{equation}\label{eqn:MD-Total-ecc}
 LE_{\varepsilon}(G)\leq n\mathrm{MD}_{\varepsilon} (G)+\E(G).
\end{equation}
Moreover the equality in \eqref{eqn:MD-Total-ecc} holds if and if $G$ is self-centred.
\end{theorem}

\begin{theorem}\label{thm:MD-Total-ecc-bi}
Let $G$ be a connected bipartite graph. Then
\begin{equation}\label{eqn:MD-Total-ecc-bi}
LE_{\varepsilon} (G)\ge \E(G).
\end{equation}
Moreover, the equality in \eqref{eqn:MD-Total-ecc-bi} holds if and only if $G$ is self-centred.
\end{theorem}

\begin{theorem}\label{thm:MD-Total-ecc-bi-low}
Let $G$ be a connected bipartite graph with $n$ vertices.  Then
\begin{equation}\label{eqn:MD-Total-ecc-bi-low}
\max\Big\{n\mathrm{MD}_{\varepsilon} (G), \E(G)\Big\}\leq LE_{\varepsilon} (G)\leq n\mathrm{MD}_{\varepsilon} (G)+\E(G).
\end{equation}
\end{theorem}

Note that a tree is a connected bipartite graph and therefore in this paper, the hypothesis "connected bipartite graph" could be replaced by "Tree".

A graph $G$ is called vertex-transitive if  for
every two vertices $u$ and $v$ of $G$, there exists an automorphism $\sigma$ of $G$ such that $\sigma( u ) = \sigma( v )$. It is known that any vertex-transitive graph is vertex degree regular, transmission regular and self-centred. Hence it follows that

\begin{corollary}\label{cor:vt}
Let $G$ be a connected vertex-transitive graph. Then the equality holds in \eqref{eqn:MD-Lap-Tr},\eqref{eqn:MD-Lap-1Zagreb},\eqref{eqn:MD-1Zagreb1},\eqref{eqn:MD-Total-ecc}. In fact,
\[\E(G)=LE_{\Tr}(G)=LE_{t}(G)=LE_{deg^2}(G)=LE_{deg}(G)=LE_{\varepsilon}(G).\]
\end{corollary}

\begin{figure}[H]
  \centering
  \includegraphics[width=5cm]{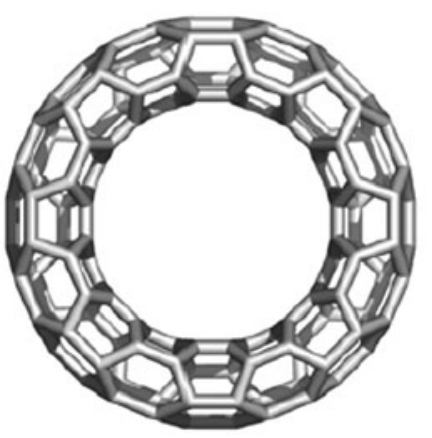}
  \caption{A toroidal fullerene (or achiral polyhex nanotorus) $T[p,q]$}\label{fig:apn}
\end{figure}
\begin{figure}[H]
  \centering
  \includegraphics[width=5cm]{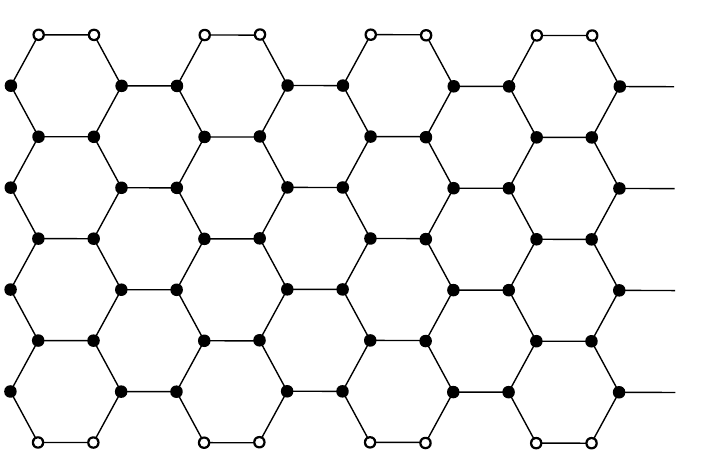}
  \caption{A 2-dimensional lattice for an achiral polyhex nanotorus $T[p,q]$}\label{fig:apnlattice}
\end{figure}

A nanostructure is an object of intermediate size between molecular and microscopic
structures. It is a product derived through engineering at the molecular
scale. In what follows we aim to apply Corollary \ref{cor:vt} to the molecular graph of a nanostructure called 
toroidal fullerenes (or achiral polyhex nanotorus) (see Fig. \ref{fig:apn} and Fig. \ref{fig:apnlattice}).

The following Lemma was proved in \cite{Ashrafi-polyhex} and \cite{Yousefi-2008}.
\begin{lemma}\label{lem:pn-vt}
The molecular graph of a polyhex nanotorus is vertex transitive.
\end{lemma}
The following is a direct consequence of Corollary \ref{cor:vt} and Lemma \ref{lem:pn-vt}.
\begin{corollary}
Let $T$ be a molecular graph of a polyhex nanotorus. Then 
\[\E(T)=LE_{\Tr}(T)=LE_{t}(T)=LE_{deg^2}(T)=LE_{deg}(T)=LE_{\varepsilon}(T).\]
\end{corollary}

\noindent\textbf{Concluding remarks:} In this paper by considering some vertex weights, some topological indices appears in Laplcian graph energy and average weight. Note that several other vertex weight and thus their corresponding Laplacian graph energy could be defined. For example if we define a weight of an arbitrary vertex $u$ as $deg(u)^3$ whose average one is $\frac{F(G)}{n}$, where $F(G)$  is  referred to as forgotten Zagreb index.

\end{document}